\documentstyle[12pt]{article}
 \normalsize

\begin{document}
\title{{\bf {On the mean value of the generalized Dirichlet $L$-functions with the weight of the Gauss Sums}}
\footnotetext{\\
1.marong@nwpu.edu.cn \qquad 2. nyn0902@mail.nwpu.edu.cn}
\author{Rong Ma$^{1}$\qquad Yana Niu$^{2}$
\\{\small{{School of Mathematics and Statistics, Northwestern
Polytechnical University}} }\\{\small{Xi'an, Shaanxi, 710072,
People's Republic of China}} \\}
\date{}}
\maketitle
\vspace{0.0cm}

\begin{center}
\large{{\bf Abstract}}
\end{center}

Let $q\ge3$ be an integer, $\chi$ denote a
Dirichlet character modulo $q$, for any real number $a\ge 0$, we
define the generalized Dirichlet $L$-functions
$$
L(s,\chi,a)=\sum_{n=1}^{\infty}\frac{\chi(n)}{(n+a)^s},
$$
where $s=\sigma+it$ with $\sigma>1$ and $t$ both real. It can be
extended to all $s$ by analytic continuation. For any integer $m$, the famous Gauss sum $G(m,\chi)$ is defined as follows:
$$G(m,\chi)=\sum_{a=1}^{q}\chi(a)e\left(\frac{am}{q}\right),
$$
where $e(y)=e^{2\pi iy}$. The main purpose of this paper is to use the analytic method to study the mean value properties of the generalized Dirichlet $L$-functions with the weight of the Gauss Sums, and obtain a sharp asymptotic formula.

\vspace{1cm}
\begin{center}
\large{{\bf1. Introduction}}
\end{center}

Let $q\ge3$ be an integer, $\chi$ denote a Dirichlet character modulo
$q$. For any interger $m$, the famous Gauss sum $G(m,\chi)$ is defined as follows:
$$G(m,\chi)=\sum_{a=1}^{q}\chi(a)e\left(\frac{am}{q}\right),
$$
where $e(y)=e^{2\pi iy}$. In particular for $m=1$, we write $\tau(\chi)=\sum_{a=1}^{q}\chi(a)e(\frac{a}{q})$.

The various properties and applications of $\tau(\chi)$ have appeared in many analytic number theory books (see Ref. [1]-[3]). Maybe the most important property of $\tau(\chi)$ is that if $\chi$ is a primitive character modulo $q$, then $|\tau(\chi)|=\sqrt q$. If $\chi$ is a non-primitive character modulo $q$, then the values distribution of $\tau(\chi)$ is more irregular, even more it may be zero! But, it is surprising that $\tau(\chi)$ appears to have many good value distribution properties in some problems of the weighted mean value. In this paper, we want to show this point.

Now let $a\ge0$ be an integer, generalized Dirichlet $L$-functions
$L(s,\chi,a)$ defined by
$$
L(s,\chi,a)=\sum_{n=1}^{\infty}\frac{\chi(n)}{(n+a)^s},
$$
where $s=\sigma+it$ with $\sigma>1$
and $t$ both real.

About the generalized Dirichlet series, B.
C. Berndt (see Ref. [4]-[6]) studied many identical properties
satisfying restrictive conditions. The first author (see Ref. [7]) also got the following asymptotic formula about the generalized Dirichlet $L$-functions
\begin{eqnarray}
\sum_{\chi\ne\chi_0}\left|L(1,\chi,a)\right|^2
=\phi(q)\sum_{d|q}\frac{\mu(d)}{d^2}\zeta\left(2,\frac{a}{d}\right)
-\frac{4\phi(q)}{a}\sum_{d|q}\frac{\mu(d)}{d}\sum_{k=1}^{[\frac{a}{d}]}\frac{1}{k}
+O\left(\frac{\phi(q)\log q}{\sqrt{q}}\right),\nonumber
\end{eqnarray}
where $\zeta(s,\alpha)(s=\sigma+it,\,\alpha>0)$ is the Hurwitz zeta function defined for $\sigma>1$ by the series
$$
\zeta(s,\alpha)=\sum_{n=0}^{\infty}\frac{1}{(n+\alpha)^{s}},
$$
and $\phi$ is the Euler function, $\mu$ is the M\"{o}bius
function, the $O$ constant only depends on $a$.

The first and second authors (see Ref. [8]) also study
the mean value properties of the generalized Dirichlet $L$-functions
with Dirichlet character and the generalized trigonometric sums.

In this paper, we will study
the mean value properties of the generalized Dirichlet $L$-functions
with the weight of the Gauss sums
$$
\sum_{\chi\ne\chi_0}|\tau(\chi)|^{m}
|L(1,\chi,a)|^{2},
$$
where $\sum_{\chi\ne\chi_0}$ denotes the summation over all non-principal characters modulo $q$, $m$ is any non-negative real number. So far, no one seems to have given any results on this content, and we don't even know if it has a good mean distribution property. In particular, we shall prove the following main conclusion:

\noindent \textbf{Theorem. }Let $q\geq3$ be an integer and $q=MN$, $M=\prod_{p||q}p$, $(M,N)=1$, $\chi$ denote a Dirichlet character modulo $q$. Then for any real number $m>0$ and positive real number $a\geq1$, we have the asymptotic formula
\begin{eqnarray}
&&\sum_{\chi\ne\chi_0}|\tau(\chi)|^{m}|L(1,\chi,a)|^{2}\nonumber\\
&=&N^{\frac{m}{2}-1}\phi^{2}(N)\zeta\left(2,\frac{a}{k}\right)
\prod_{p||q}(1-\frac{1}{p^2})\prod_{p|M}\left(p^{\frac{m}{2}+1}
-2p^{\frac{m}{2}}+1\right)+O\left(q^{\frac{m}{2}+\epsilon}\right),\nonumber
\end{eqnarray}
where $\prod_{p||q}$ denotes the production over all prime divisors satisfied with $p|q$ and $p^{2}\dag q$, $\phi$ is the Euler function,
$\zeta(s,\alpha)(s=\sigma+it,\,\alpha>0)$ is the Hurwitz zeta
function, $\epsilon$ is any fixed positive number.

\begin{center}
\large{{\bf 2. Some lemmas}}
\end{center}

To complete the proof of the Theorem, we need the following several
lemmas. First, we make an identity of
the Dirichlet $L$-functions and the generalized form.

\noindent \textbf{Lemma 1. }Let $q\ge3$ be an integer, and $\chi$
denote a non-principal Dirichlet character modulo $q$. Let $L(s,\chi)$ denote the
Dirichlet $L$-functions corresponding to $\chi$, and $L(s,\chi,a)$
denote the generalized Dirichlet $L$-functions. Then for any real
number $a\ge0$, we have
$$
L(1,\chi,a)=L(1,\chi)-a\sum_{n=1}^{\infty}\frac{\chi(n)}{n(n+a)}.
$$
\\
\textbf{Proof. }See Lemma 1, and let $m=1$ (Ref. [7]).

\noindent \textbf{Lemma 2. }Let $q=uv$, $u\geq2$, $v\geq2$, $(u,v)=1$. Then for any character $\chi$ modulo $q$, there exist a unique character $\chi_u$ modulo $u$ and a unique character $\chi_v$ modulo $v$ such that $\chi=\chi_u\chi_v$, and
$$
|\tau(\chi)|=|\tau(\chi_u)|\times|\tau(\chi_v)|.
$$
\\
\textbf{Proof. }See Theorem 13.3.1 of reference [2].

\noindent \textbf{Lemma 3. }Let $p$ be a prime, $\alpha$ be a positive integer with $\alpha\geq2$, $n=p^{\alpha}$. Then for any non-primitive character $\chi_1$ modulo $n$, we have the identity
$$
\tau(\chi_1)=\sum_{a=1}^{p^{\alpha}}\chi_1(a)e\left(\frac{a}{p^{\alpha}}\right)=0.
$$
\\
\textbf{Proof. }For any non-primitive character $\chi_1$
modulo $n=p^{\alpha}$, it is clear that $\chi_1$ is also a character modulo $p^{\alpha-1}$. if b and l round through a complete residual system modulo $p^{\alpha-1}$ and modulo $p$, respectively, then $b+p^{\alpha-1}l$ also round through a complete residual system modulo $p^{\alpha}$; so noting that $\chi_1$ is a character modulo $p^{\alpha-1}$ and the identity $\sum_{l=0}^{p-1}e\left(\frac{l}{p}\right)=0$, we have
\begin{eqnarray}
&&\sum_{a=1}^{p^{\alpha}}\chi_1(a)e\left(\frac{a}{p^{\alpha}}\right)\nonumber\\
&=&\sum_{l=0}^{p-1}\sum_{b=1}^{p^{\alpha-1}}\chi_1(b+p^{\alpha-1}l)
e\left(\frac{b+p^{\alpha-1}l}{p^{\alpha}}\right)\nonumber\\
&=&\left(\sum_{b=1}^{p^{\alpha-1}}\chi_1(b)e\left(\frac{b}{p^{\alpha}}\right)\right)
\left(\sum_{l=0}^{p-1}e\left(\frac{l}{p}\right)\right)\nonumber\\
&=&0.\nonumber
\end{eqnarray}
This completes the proof of Lemma 3.

\noindent \textbf{Lemma 4. }Let $q$ and $r$ be integers with $q\ge2$ and $(r,q)=1$, and $\chi$ be a Dirichlet character modulo $q$. Then we have the identities
$$
\displaystyle\mathop{{\sum}^{\ast}}_{\chi \bmod q}\chi(r)=\sum_{d|(q,r-1)}\mu\left(\frac{q}{d}\right)\phi(d) \quad and  \quad
J(q)=\sum_{d|q}\mu(d)\phi\left(\frac{q}{d}\right),
$$
where $\mu(n)$ is a M\"{o}bius function, $\phi(n)$ is a Euler function, and $J(q)$ denotes the number of all primitive characters modulo $q$.\\
\\
\textbf{Proof. }See lemma 3 of reference [9].

\noindent \textbf{Lemma 5. }Let $q\geq3$ be an integer and $q=MN$, $M=\prod_{p||q}p$, $(M,N)=1$, $\chi$ denote a Dirichlet character modulo $q$. Then for any real number $m>0$, we have the asymptotic formula
\begin{eqnarray}
&&\sum_{\chi\ne\chi_0}|\tau(\chi)|^{m}|L(1,\chi)|^{2}\nonumber\\
&=&N^{\frac{m}{2}-1}\phi^{2}(N)\zeta(2)\prod_{p||q}(1-\frac{1}{p^2})\prod_{p|M}\left(p^{\frac{m}{2}+1}
-2p^{\frac{m}{2}}+1\right)+O\left(q^{\frac{m}{2}+\epsilon}\right),\nonumber
\end{eqnarray}
where $\prod_{p||q}$ denotes the production over all prime divisors satisfied with $p|q$ and $p^{2}\dag q$, $\phi$ is the Euler function, $\epsilon$ is any fixed positive number.\\
\\
\textbf{Proof. }See Theorem, and let $k=1$ (Ref. [10]).

\noindent \textbf{Lemma 6. }Let $q=p_1\dots p_kp_{k+1}^{\alpha_{k+1}}\dots p_r^{\alpha_r}$, where $\alpha_i>1$, $k+1\leq i\leq r$. Let $M=p_1p_2\dots p_k$, $N=p_{k+1}^{\alpha_{k+1}}\dots p_r^{\alpha_r}$, then $(M,N)=1$. For any positive integer $n$ and $d|M$, we have
\begin{eqnarray}
&&\displaystyle\mathop{{\sum}^{\ast}}_{\chi \bmod Nd}
\sum_{n=1}^{\infty}\frac{\chi{\chi_{M}^0}(n)}{n(n+a)}L(1,\overline{\chi{\chi_{M}^0}})\nonumber\\
&=&\frac{\phi^{2}(N)}{N}J(d)\zeta(2)\sum_{k|q}\frac{\mu(k)}{ak^2}-
\frac{\phi^{2}(N)}{N}J(d)\sum_{k|q}\frac{\mu(k)}{a^2k}\sum_{h=1}^{[a/k]}\frac{1}{h}+O((Nd)^{\epsilon}),\nonumber
\end{eqnarray}
where $\sum_{\chi \bmod Nd}^{\ast}$ denotes the summation over all primitive characters modulo $Nd$, $J(d)$ denotes the number of all primitive characters modulo $d$, $\epsilon$ is any fixed positive number, $\chi_{M}^0$ denotes the principal character of modulo $M$.\\
\\
\textbf{Proof. }For convenience, we put
$$A(\chi,y)=\sum_{Nd\le n\le y}\chi(n).
$$
Then applying Abel's identity, by analytic continuation we have
$$
L(1,\bar{\chi})=\sum_{n=1}^{Nd}\frac{\bar{\chi}(n)}{n}+\int_{Nd}^{+\infty}\frac{A(\bar{\chi},y)}{y^2}dy,
$$
$$
\sum_{n=1}^{\infty}\frac{\chi(n)}{n(n+a)}=\sum_{n=1}^{Nd}\frac{\chi(n)}{n(n+a)}
+\int_{N}^{+\infty}\frac{(2y+a)A(\chi,y)}{y^2(y+a)^2}dy.
$$
So,
\begin{eqnarray}
& &\displaystyle\mathop{{\sum}^{\ast}}_{\chi \bmod Nd}
\sum_{n=1}^{\infty}\frac{\chi{\chi_{M}^0}(n)}{n(n+a)}L(1,\overline{\chi{\chi_{M}^0}})\nonumber\\
&=&\displaystyle\mathop{{\sum}^{\ast}}_{\chi \bmod Nd}
\left(\sum_{n=1}^{Nd}\frac{\chi{\chi_{M}^0}(n)}{n(n+a)}+\int_{Nd}^{+\infty}
\frac{(2y+a)A(\chi{\chi_{M}^0},y)}{y^2(y+a)^2}dy\right)\times\nonumber\\
& &\times\left(\sum_{n=1}^{Nd}\frac{\overline{\chi{\chi_{M}^0}}(n)}{n}
+\int_{Nd}^{+\infty}\frac{A(\overline{\chi{\chi_{M}^0}},y)}{y^2}dy\right)\nonumber\\
&=&\displaystyle\mathop{{\sum}^{\ast}}_{\chi \bmod Nd}
\sum_{n=1}^{Nd}\frac{\chi{\chi_{M}^0}(n)}{n(n+a)}
\sum_{m=1}^{Nd}\frac{\overline{\chi{\chi_{M}^0}}(m)}{m}+\nonumber\\
& &+\displaystyle\mathop{{\sum}^{\ast}}_{\chi \bmod Nd}
\sum_{n=1}^{Nd}\frac{\chi{\chi_{M}^0}(n)}{n(n+a)}
\int_{Nd}^{+\infty}\frac{(A(\overline{\chi{\chi_{M}^0}},y)}{y^2}dy+\nonumber\\
& &+\displaystyle\mathop{{\sum}^{\ast}}_{\chi \bmod Nd}
\sum_{n=1}^{Nd}\frac{\overline{\chi{\chi_{M}^0}}(n)}{n}
\int_{Nd}^{+\infty}\frac{(2y+a)A(\chi{\chi_{M}^0},y)}{y^2(y+a)^2}dy+\nonumber\\
& &+\displaystyle\mathop{{\sum}^{\ast}}_{\chi \bmod Nd}
\int_{Nd}^{+\infty}\frac{(2y+a)A(\chi{\chi_{M}^0},y)}{y^2(y+a)^2}dy
\int_{Nd}^{+\infty}\frac{(A(\overline{\chi{\chi_{M}^0}},z)}{z^2}dz\nonumber\\
&=&A_1+A_2+A_3+A_4.
\end{eqnarray}

We will estimate each of them. Firstly we estimate $A_1$, According to Lemma 4, we have
\begin{eqnarray}
A_1&=&\displaystyle\mathop{{\sum}^{\ast}}_{\chi \bmod Nd}
\sum_{n=1}^{Nd}\frac{\chi{\chi_{M}^0}(n)}{n(n+a)}
\sum_{m=1}^{Nd}\frac{\overline{\chi{\chi_{M}^0}}(m)}{m}\nonumber\\
&=&\sum_{n=1}^{Nd}\sum_{m=1}^{Nd}\frac{{\chi_{M}^0}(n\bar{m})}{mn(n+a)}
\displaystyle\mathop{{\sum}^{\ast}}_{\chi \bmod Nd}
\chi(n\bar{m})\nonumber\\
&=&\displaystyle\mathop{\sum{'}}_{n=1}^{Nd}
\displaystyle\mathop{\sum{'}}_{m=1}^{Nd}
\frac{{\chi_{M}^0}(n\bar{m})}{mn(n+a)}
\sum_{l|(Nd,n\bar{m}-1)}\mu(\frac{Nd}{l})\phi(l)\nonumber\\
&=&\sum_{l|Nd}\mu(\frac{Nd}{l})\phi(l)\displaystyle\mathop
{\displaystyle\mathop{\sum{'}}_{n=1}^{Nd}
\displaystyle\mathop{\sum{'}}_{m=1}^{Nd}}
_{n\equiv m(\bmod l)}\frac{{\chi_{M}^0}(n\bar{m})}{mn(n+a)}\nonumber\\
&=&\sum_{l|Nd}\mu(\frac{Nd}{l})\phi(l)
\displaystyle\mathop{\sum{'}}_{n=1}^{Nd}
\frac{1}{n^2(n+a)}+\sum_{l|Nd}\mu(\frac{Nd}{l})\phi(l)\displaystyle\mathop
{\displaystyle\mathop{\sum{'}}_{n=1}^{Nd}
\displaystyle\mathop{\sum{'}}_{m=1}^{Nd}}
_{n\equiv m(\bmod l),\,n\neq m}
\frac{{\chi_{M}^0}(n\bar{m})}{mn(n+a)}\nonumber\\
&=&\sum_{l|Nd}\mu(\frac{Nd}{l})\phi(l)\sum_{n=1}^{Nd}
\frac{1}{n^2(n+a)}\sum_{k|(n,q)}\mu(k)+\nonumber\\
&&+O\left(\sum_{l|Nd}\phi(l)
\displaystyle\mathop
{\displaystyle\mathop{\sum{'}}_{n=1}^{Nd}
\displaystyle\mathop{\sum{'}}_{m=1}^{Nd}}
_{n\equiv m(\bmod l),\,n\neq m}
\frac{1}{mn(n+a)}\right)\nonumber\\
&=&\sum_{l|Nd}\mu(\frac{Nd}{l})\phi(l)
\sum_{k|q}\mu(k)\sum_{n=1 \atop k|n}^{Nd}\frac{1}{n^2(n+a)}
+O((Nd)^{\epsilon})\nonumber\\
&=&J(Nd)\sum_{k|q}\frac{\mu(k)}{k^3}\sum_{n=1}^{Nd/k}\frac{1}{n^2(n+a/k)}
+O((Nd)^{\epsilon})\nonumber\\
&=&J(Nd)\sum_{k|q}\frac{\mu(k)}{k^3}\sum_{n=1}^{\infty}\frac{1}{n^2(n+a/k)}+\nonumber\\
& &+O\left(J(Nd)\sum_{k|q}\frac{\mu(k)}{k^3}
\sum_{n>{Nd/k}}\frac{1}{n^2(n+a/k)}\right)+O((Nd)^{\epsilon})\nonumber\\
&=&J(Nd)\sum_{k|q}\frac{\mu(k)}{k^3}\left(\sum_{n=1}^{\infty}\frac{1}{(a/k)^2}
\left(\frac{a/k}{n^2}+\frac{1}{n+a/k}-\frac{1}{n}\right)\right)
+O((Nd)^{\epsilon})\nonumber\\
&=&J(Nd)\sum_{k|q}\frac{\mu(k)}{ak^2}\sum_{n=1}^{\infty}\frac{1}{n^2}
+J(Nd)\sum_{k|q}\frac{\mu(k)}{a^2k}\sum_{n=1}^{\infty}
\left(\frac{1}{n+a/k}-\frac{1}{n}\right)+O((Nd)^{\epsilon})\nonumber\\
&=&J(Nd)\sum_{k|q}\frac{\mu(k)}{ak^2}\zeta(2)-
J(Nd)\sum_{k|q}\frac{\mu(k)}{a^2k}\sum_{h=1}^{[a/k]}\frac{1}{h}+O((Nd)^{\epsilon})\nonumber\\
&=&\frac{\phi^{2}(N)}{N}J(d)\zeta(2)\sum_{k|q}\frac{\mu(k)}{ak^2}-
\frac{\phi^{2}(N)}{N}J(d)\sum_{k|q}\frac{\mu(k)}{a^2k}\sum_{h=1}^{[a/k]}\frac{1}{h}+O((Nd)^{\epsilon}),
\end{eqnarray}
where $\sum{'}_{n}$ indicates that the sum is over those $n$ relatively prime to $q$.

In the following we will estimate $A_2,\,A_3$ and $A_4$, According to Cauchy inequality and Polya-Vinogradiv inequality about character sums we can easily get
\begin{eqnarray}
A_2&=&\displaystyle\mathop{{\sum}^{\ast}}_{\chi \bmod Nd}
\sum_{n=1}^{Nd}\frac{\chi{\chi_{M}^0}(n)}{n(n+a)}
\int_{Nd}^{+\infty}\frac{A(\overline{\chi{\chi_{M}^0}},y)}{y^2}dy\nonumber\\
&=&\displaystyle\mathop{{\sum}^{\ast}}_{\chi \bmod Nd}
\sum_{n=1}^{Nd}\frac{\chi(n){\chi_{M}^0}(n)}{n(n+a)}
\int_{Nd}^{(Nd)^{3\cdot2^{n-2}}}\frac{\sum_{Nd<m\leq y}\bar{\chi}(m)\overline{\chi_{M}^0}(m)}{y^2}dy+\nonumber\\
& &+\displaystyle\mathop{{\sum}^{\ast}}_{\chi \bmod Nd}
\sum_{n=1}^{Nd}\frac{\chi(n){\chi_{M}^0}(n)}{n(n+a)}
\int_{(Nd)^{3\cdot2^{n-2}}}^{\infty}\frac{\sum_{Nd<m\leq y}\bar{\chi}(m)\overline{\chi_{M}^0}(m)}{y^2}dy\nonumber\\
&\ll&\int_{Nd}^{(Nd)^{3\cdot2^{n-2}}}\frac{1}{y^2}\left|
\sum_{n=1}^{Nd}\sum_{m=Nd}^{y}\frac{1}{n(n+a)}
\displaystyle\mathop{{\sum}^{\ast}}_{\chi \bmod Nd}\chi(n\bar{m})\right|dy+\nonumber\\
&&+(Nd)^{\epsilon}\int_{(Nd)^{3\cdot2^{n-2}}}^{\infty}\frac{1}{y^2}
\displaystyle\mathop{{\sum}^{\ast}}_{\chi \bmod Nd}
|A(\bar{\chi},y)|dy\nonumber\\
&\leq&\int_{Nd}^{(Nd)^{3\cdot2^{n-2}}}\frac{1}{y^2}\left|
\sum_{n=1}^{Nd}\sum_{m=Nd}^{y}\frac{1}{n(n+a)}
\sum_{l|(Nd,n\bar{m}-1)}\phi(l)\right|dy+\nonumber\\
&&+O\left((Nd)^{\frac{3}{2}+\epsilon}
\int_{(Nd)^{3\cdot2^{n-2}}}^{\infty}\frac{1}{y^2}dy\right)\nonumber\\
&\leq&\int_{Nd}^{(Nd)^{3\cdot2^{n-2}}}\frac{1}{y^2}\left|
\sum_{l|Nd}\phi(l)\displaystyle\mathop
{\displaystyle\mathop{\sum{'}}_{n=1}^{Nd}
\displaystyle\mathop{\sum{'}}_{m=Nd}^{y}}
_{n\equiv m(\bmod l)}\frac{1}{n(n+a)}\right|dy
+O((Nd)^{\epsilon})\nonumber\\
&\leq&\int_{Nd}^{(Nd)^{3\cdot2^{n-2}}}\frac{1}{y^2}\left|
\sum_{l|Nd}\phi(l)\displaystyle\mathop{\sum{'}}_{n=1}^{Nd}\frac{1}{n(n+a)}
\cdot\frac{y}{l}\cdot(Nd)^{\epsilon}\right|dy+O((Nd)^{\epsilon})\nonumber\\
&\ll&(Nd)^{\epsilon},
\end{eqnarray}
Similarly, we also have
\begin{eqnarray}
A_3&=&\displaystyle\mathop{{\sum}^{\ast}}_{\chi \bmod Nd}
\sum_{n=1}^{Nd}\frac{\overline{\chi{\chi_{M}^0}}(n)}{n}
\int_{Nd}^{+\infty}\frac{(2y+a)A(\chi{\chi_{M}^0},y)}{y^2(y+a)^2}dy\ll(Nd)^{\epsilon},
\end{eqnarray}
\begin{eqnarray}
A_4&=&\displaystyle\mathop{{\sum}^{\ast}}_{\chi \bmod Nd}
\int_{Nd}^{+\infty}\frac{(2y+a)A(\chi{\chi_{M}^0},y)}{y^2(y+a)^2}dy
\int_{Nd}^{+\infty}\frac{(A(\overline{\chi{\chi_{M}^0}},z)}{z^2}dz\ll(Nd)^{\epsilon}.
\end{eqnarray}

Combining the formulas (1)-(5), we immediately obtain
\begin{eqnarray}
&&\displaystyle\mathop{{\sum}^{\ast}}_{\chi \bmod Nd}
\sum_{n=1}^{\infty}\frac{\chi{\chi_{M}^0}(n)}{n(n+a)}L(1,\overline{\chi{\chi_{M}^0}})\nonumber\\
&=&\frac{\phi^{2}(N)}{N}J(d)\zeta(2)\sum_{k|q}\frac{\mu(k)}{ak^2}-
\frac{\phi^{2}(N)}{N}J(d)\sum_{k|q}\frac{\mu(k)}{a^2k}\sum_{h=1}^{[a/k]}\frac{1}{h}+O((Nd)^{\epsilon}).\nonumber
\end{eqnarray}
This completes the proof of Lemma 6.

\noindent \textbf{Lemma 7. }Let $q=p_1\dots p_kp_{k+1}^{\alpha_{k+1}}\dots p_r^{\alpha_r}$, where $\alpha_i>1$, $k+1\leq i\leq r$. Let $M=p_1p_2\dots p_k$, $N=p_{k+1}^{\alpha_{k+1}}\dots p_r^{\alpha_r}$, then $(M,N)=1$. For any positive integer $n$ and $d|M$, we have
\begin{eqnarray}
&&\displaystyle\mathop{{\sum}^{\ast}}_{\chi \bmod Nd}
\left|\sum_{n=1}^{\infty}\frac{\chi{\chi_{M}^0}(n)}{n(n+a)}\right|^2\nonumber\\
&=&\frac{\phi^{2}(N)}{N}J(d)\zeta(2)\sum_{k|q}\frac{\mu(k)}{a^2k^2}
+\frac{\phi^{2}(N)}{N}J(d)\zeta(2,\frac{a}{k})\sum_{k|q}\frac{\mu(k)}{a^2k^2}-\nonumber\\
&&-2\frac{\phi^{2}(N)}{N}J(d)\sum_{k|q}\frac{\mu(k)}{a^3k}
\sum_{h=1}^{[a/k]}\frac{1}{h}+O((Nd)^{\epsilon}),\nonumber
\end{eqnarray}
where $\sum_{\chi \bmod Nd}^{\ast}$ denotes the summation over all primitive characters modulo $Nd$, $J(d)$ denotes the number of all primitive characters modulo $d$, $\epsilon$ is any fixed positive number, $\chi_{M}^0$ denotes the principal character of modulo $M$.\\
\\
\textbf{Proof. }According to Abel's identity, and $A(\chi,y)=\sum_{Nd\le n\le y}\chi(n)$ as defined in the proof of Lemma 6, we have
\begin{eqnarray}
&&\displaystyle\mathop{{\sum}^{\ast}}_{\chi \bmod Nd}
\left|\sum_{n=1}^{\infty}\frac{\chi{\chi_{M}^0}(n)}{n(n+a)}\right|^2\nonumber\\
&=&\displaystyle\mathop{{\sum}^{\ast}}_{\chi \bmod Nd}
\left(\sum_{n=1}^{Nd}\frac{\chi{\chi_{M}^0}(n)}{n(n+a)}
+\int_{Nd}^{+\infty}\frac{(2y+a)A(\chi\chi_{M}^0,y)}{y^2(y+a)^2}dy\right)\times\nonumber\\
&&\times\left(\sum_{m=1}^{Nd}\frac{\overline{\chi\chi_{M}^0}(m)}{m(m+a)}
+\int_{Nd}^{+\infty}\frac{(2z+a)A(\overline{\chi\chi_{M}^0},z)}{z^2(z+a)^2}dz\right)\nonumber\\
&=&\displaystyle\mathop{{\sum}^{\ast}}_{\chi \bmod Nd}
\left(\sum_{n=1}^{Nd}\frac{\chi{\chi_{M}^0}(n)}{n(n+a)}\right)
\left(\sum_{m=1}^{Nd}\frac{\overline{\chi\chi_{M}^0}(m)}{m(m+a)}\right)+\nonumber\\
&&+\displaystyle\mathop{{\sum}^{\ast}}_{\chi \bmod Nd}
\left(\sum_{n=1}^{Nd}\frac{\chi{\chi_{M}^0}(n)}{n(n+a)}\right)
\left(\int_{Nd}^{+\infty}\frac{(2z+a)A(\overline{\chi\chi_{M}^0},z)}{z^2(z+a)^2}dz\right)+\nonumber\\
&&+\displaystyle\mathop{{\sum}^{\ast}}_{\chi \bmod Nd}
\left(\int_{Nd}^{+\infty}\frac{(2y+a)A(\chi\chi_{M}^0,y)}{y^2(y+a)^2}dy\right)
\left(\sum_{m=1}^{Nd}\frac{\overline{\chi\chi_{M}^0}(m)}{m(m+a)}\right)+\nonumber\\
&&+\displaystyle\mathop{{\sum}^{\ast}}_{\chi \bmod Nd}
\left(\int_{Nd}^{+\infty}\frac{(2y+a)A(\chi\chi_{M}^0,y)}{y^2(y+a)^2}dy\right)
\left(\int_{Nd}^{+\infty}\frac{(2z+a)A(\overline{\chi\chi_{M}^0},z)}{z^2(z+a)^2}dz\right)\nonumber\\
&=&B_1+B_2+B_3+B_4.
\end{eqnarray}

We will estimate each of them. Firstly we estimate $B_1$. According to Lemma 4, we have
\begin{eqnarray}
B_1&=&\displaystyle\mathop{{\sum}^{\ast}}_{\chi \bmod Nd}
\left(\sum_{n=1}^{Nd}\frac{\chi{\chi_{M}^0}(n)}{n(n+a)}\right)
\left(\sum_{m=1}^{Nd}\frac{\overline{\chi\chi_{M}^0}(m)}{m(m+a)}\right)\nonumber\\
&=&\sum_{n=1}^{Nd}\sum_{m=1}^{Nd}\frac{{\chi_{M}^0}(n\bar{m})}{mn(m+a)(n+a)}
\displaystyle\mathop{{\sum}^{\ast}}_{\chi \bmod Nd}
\chi(n\bar{m})\nonumber\\
&=&\displaystyle\mathop{\sum{'}}_{n=1}^{Nd}
\displaystyle\mathop{\sum{'}}_{m=1}^{Nd}
\frac{{\chi_{M}^0}(n\bar{m})}{mn(m+a)(n+a)}
\sum_{l|(Nd,n\bar{m}-1)}\mu(\frac{Nd}{l})\phi(l)\nonumber\\
&=&\sum_{l|Nd}\mu(\frac{Nd}{l})\phi(l)\displaystyle\mathop
{\displaystyle\mathop{\sum{'}}_{n=1}^{Nd}
\displaystyle\mathop{\sum{'}}_{m=1}^{Nd}}
_{n\equiv m(\bmod l)}\frac{{\chi_{M}^0}(n\bar{m})}{mn(m+a)(n+a)}\nonumber\\
&=&\sum_{l|Nd}\mu(\frac{Nd}{l})\phi(l)
\displaystyle\mathop{\sum{'}}_{n=1}^{Nd}
\frac{1}{n^2(n+a)^2}+\nonumber\\
& &+\sum_{l|Nd}\mu(\frac{Nd}{l})\phi(l)\displaystyle\mathop
{\displaystyle\mathop{\sum{'}}_{n=1}^{Nd}
\displaystyle\mathop{\sum{'}}_{m=1}^{Nd}}
_{n\equiv m(\bmod l),\,n\neq m}
\frac{{\chi_{M}^0}(n\bar{m})}{mn(m+a)(n+a)}\nonumber\\
&=&\sum_{l|Nd}\mu(\frac{Nd}{l})\phi(l)\sum_{n=1}^{Nd}
\frac{1}{n^2(n+a)^2}\sum_{k|(n,q)}\mu(k)+\nonumber\\
& &+O\left(\sum_{l|Nd}\phi(l)
\displaystyle\mathop
{\displaystyle\mathop{\sum{'}}_{n=1}^{Nd}
\displaystyle\mathop{\sum{'}}_{m=1}^{Nd}}
_{n\equiv m(\bmod l),\,n\neq m}
\frac{1}{mn(m+a)(n+a)}\right)\nonumber\\
&=&\sum_{l|Nd}\mu(\frac{Nd}{l})\phi(l)
\sum_{k|q}\mu(k)\sum_{n=1 \atop k|n}^{Nd}\frac{1}{n^2(n+a)^2}
+O((Nd)^{\epsilon})\nonumber\\
&=&J(Nd)\sum_{k|q}\frac{\mu(k)}{k^4}\sum_{n=1}^{Nd/k}\frac{1}{n^2(n+a/k)^2}
+O((Nd)^{\epsilon})\nonumber\\
&=&J(Nd)\sum_{k|q}\frac{\mu(k)}{k^4}\sum_{n=1}^{\infty}\frac{1}{n^2(n+a/k)^2}+\nonumber\\
& &+O\left(J(Nd)\sum_{k|q}\frac{\mu(k)}{k^4}
\sum_{n>{Nd/k}}\frac{1}{n^2(n+a/k)^2}\right)+O((Nd)^{\epsilon})\nonumber\\
&=&J(Nd)\sum_{k|q}\frac{\mu(k)}{k^4}\left(\sum_{n=1}^{\infty}\frac{1}{(a/k)^2}
\left(\frac{1}{n^2}+\frac{1}{(n+a/k)^2}-\frac{2k/a}{n}+\frac{2k/a}{n+a/k}\right)\right)+\nonumber\\
& &+O((Nd)^{\epsilon})\nonumber\\
&=&J(Nd)\sum_{k|q}\frac{\mu(k)}{a^2k^2}\sum_{n=1}^{\infty}\frac{1}{n^2}
+J(Nd)\sum_{k|q}\frac{\mu(k)}{a^2k^2}\sum_{n=1}^{\infty}\frac{1}{(n+a/k)^2}-\nonumber\\
& &-2J(Nd)\sum_{k|q}\frac{\mu(k)}{a^3k}\sum_{n=1}^{\infty}
\left(\frac{1}{n}-\frac{1}{n+a/k}\right)+O((Nd)^{\epsilon})\nonumber\\
&=&J(Nd)\zeta(2)\sum_{k|q}\frac{\mu(k)}{a^2k^2}+J(Nd)\zeta(2,\frac{a}{k})\sum_{k|q}\frac{\mu(k)}{a^2k^2}-\nonumber\\
&&-2J(Nd)\sum_{k|q}\frac{\mu(k)}{a^3k}
\sum_{h=1}^{[a/k]}\frac{1}{h}+O((Nd)^{\epsilon})\nonumber\\
&=&\frac{\phi^{2}(N)}{N}J(d)\zeta(2)\sum_{k|q}\frac{\mu(k)}{a^2k^2}
+\frac{\phi^{2}(N)}{N}J(d)\zeta(2,\frac{a}{k})\sum_{k|q}\frac{\mu(k)}{a^2k^2}-\nonumber\\
&&-2\frac{\phi^{2}(N)}{N}J(d)\sum_{k|q}\frac{\mu(k)}{a^3k}
\sum_{h=1}^{[a/k]}\frac{1}{h}+O((Nd)^{\epsilon}),
\end{eqnarray}
we also could get
\begin{eqnarray}
B_2&=&\displaystyle\mathop{{\sum}^{\ast}}_{\chi \bmod Nd}
\left(\sum_{n=1}^{Nd}\frac{\chi{\chi_{M}^0}(n)}{n(n+a)}\right)
\left(\int_{Nd}^{+\infty}\frac{(2z+a)A(\overline{\chi\chi_{M}^0},z)}{z^2(z+a)^2}dz\right)\ll(Nd)^{\epsilon},
\end{eqnarray}
\begin{eqnarray}
B_3&=&\displaystyle\mathop{{\sum}^{\ast}}_{\chi \bmod Nd}
\left(\int_{Nd}^{+\infty}\frac{(2y+a)A(\chi\chi_{M}^0,y)}{y^2(y+a)^2}dy\right)
\left(\sum_{m=1}^{Nd}\frac{\overline{\chi\chi_{M}^0}(m)}{m(m+a)}\right)\ll(Nd)^{\epsilon},
\end{eqnarray}
\begin{eqnarray}
B_4&=&\displaystyle\mathop{{\sum}^{\ast}}_{\chi \bmod Nd}
\left(\int_{Nd}^{+\infty}\frac{(2y+a)A(\chi\chi_{M}^0,y)}{y^2(y+a)^2}dy\right)
\left(\int_{Nd}^{+\infty}\frac{(2z+a)A(\overline{\chi\chi_{M}^0},z)}{z^2(z+a)^2}dz\right)\nonumber\\
&\ll&(Nd)^{\epsilon}.
\end{eqnarray}

Combining the formulas (6)-(10), we immediately obtain
\begin{eqnarray}
&&\displaystyle\mathop{{\sum}^{\ast}}_{\chi \bmod Nd}
\left|\sum_{n=1}^{\infty}\frac{\chi{\chi_{M}^0}(n)}{n(n+a)}\right|^2\nonumber\\
&=&\frac{\phi^{2}(N)}{N}J(d)\zeta(2)\sum_{k|q}\frac{\mu(k)}{a^2k^2}
+\frac{\phi^{2}(N)}{N}J(d)\zeta(2,\frac{a}{k})\sum_{k|q}\frac{\mu(k)}{a^2k^2}-\nonumber\\
&&-2\frac{\phi^{2}(N)}{N}J(d)\sum_{k|q}\frac{\mu(k)}{a^3k}
\sum_{h=1}^{[a/k]}\frac{1}{h}+O((Nd)^{\epsilon}).\nonumber
\end{eqnarray}
This completes the proof of Lemma 7.

\begin{center}
\large{{\bf 3. Proof of Theorem}}
\end{center}

In this section, we shall complete the proof of the theorem. First, from Lemma 1, we get
\begin{eqnarray}
&&\sum_{\chi\ne\chi_0}|\tau(\chi)|^{m}|L(1,\chi,a)|^{2}\nonumber\\
&=&\sum_{\chi\ne\chi_0}|\tau(\chi)|^{m}\left|L(1,\chi)-a\sum_{n=1}^{\infty}\frac{\chi(n)}{n(n+a)}\right|^2 \nonumber\\
&=&\sum_{\chi\ne\chi_0}|\tau(\chi)|^{m}|L(1,\chi)|^{2}-a\sum_{\chi\ne\chi_0}|\tau(\chi)|^{m}
\sum_{n=1}^{\infty}\frac{\chi(n)}{n(n+a)}L(1,\bar{\chi})-\nonumber\\
& &-a\sum_{\chi\ne\chi_0}|\tau(\chi)|^{m}
\sum_{n=1}^{\infty}\frac{\bar{\chi}(n)}{n(n+a)}L(1,\chi)+a^2\sum_{\chi\ne\chi_0}|\tau(\chi)|^{m}
\left|\sum_{n=1}^{\infty}\frac{\chi(n)}{n(n+a)}\right|^2\nonumber\\
&=&C_1-aC_2-aC_3+a^2C_4,\nonumber
\end{eqnarray}
where
\begin{eqnarray}
&C_1&=\sum_{\chi\ne\chi_0}|\tau(\chi)|^{m}|L(1,\chi)|^{2},\nonumber\\
&C_2&=\sum_{\chi\ne\chi_0}|\tau(\chi)|^{m}
\sum_{n=1}^{\infty}\frac{\chi(n)}{n(n+a)}L(1,\bar{\chi}),\nonumber\\
&C_3&=\sum_{\chi\ne\chi_0}|\tau(\chi)|^{m}
\sum_{n=1}^{\infty}\frac{\bar{\chi}(n)}{n(n+a)}L(1,\chi),\nonumber\\
&C_4&=\sum_{\chi\ne\chi_0}|\tau(\chi)|^{m}
\left|\sum_{n=1}^{\infty}\frac{\chi(n)}{n(n+a)}\right|^2.\nonumber
\end{eqnarray}
From Lemma 5, we could get $C_1$ immediately. Therefore, we will estimate $C_2$, $C_3$ and $C_4$ later.

(i) The asymptotic formula of $C_2$

Let $q=p_1\dots p_kp_{k+1}^{\alpha_{k+1}}\dots p_r^{\alpha_r}$, where $\alpha_i>1$, $k+1\leq i\leq r$. Let $M=p_1p_2\dots p_k$, $N=p_{k+1}^{\alpha_{k+1}}\dots p_r^{\alpha_r}$, then $(M,N)=1$.
From the properties of primitive characters, Lemma 2 and Lemma 3, we have
\begin{eqnarray}
C_2&=&\sum_{\chi\ne\chi_0}|\tau(\chi)|^{m}\sum_{n=1}^{\infty}\frac{\chi(n)}{n(n+a)}L(1,\bar{\chi})\nonumber\\
&=&\displaystyle\mathop{\displaystyle\mathop\sum_{\chi_1\bmod p_1}\dots\displaystyle\mathop\sum_{\chi_k\bmod p_k}\displaystyle\mathop\sum_{\chi_{k+1}\bmod p_{k+1}^{\alpha_{k+1}}}\dots
\displaystyle\mathop\sum_{\chi_r\bmod p_r^{\alpha_r}}}_{\atop \chi_1\dots\chi_k\chi_{k+1}\dots\chi_r\neq \chi_0}\left|\tau(\chi_1\dots\chi_k\chi_{k+1}\dots \chi_r)\right|^m\times\nonumber\\
& &\times\sum_{n=1}^{\infty}\frac{\chi_1\dots \chi_r(n)}{n(n+a)}L(1,\overline{\chi_1\dots\chi_r})\nonumber\\
&=&\displaystyle\mathop{\displaystyle\mathop\sum_{\chi_1\bmod p_1}\dots\displaystyle\mathop\sum_{\chi_k\bmod p_k}\displaystyle\mathop\sum_{\chi_{k+1}\bmod p_{k+1}^{\alpha_{k+1}}}\dots
\displaystyle\mathop\sum_{\chi_r\bmod p_r^{\alpha_r}}}_{\atop \chi_1\dots\chi_k\chi_{k+1}\dots\chi_r\neq \chi_0}|\tau(\chi_1)|^m\dots|\tau(\chi_k)|^m\cdot\nonumber\\
&&\cdot|\tau(\chi_{k+1})|^m\dots|\tau(\chi_r)|^m
\sum_{n=1}^{\infty}\frac{\chi_1\dots
\chi_r(n)}{n(n+a)}L(1,\overline{\chi_1\dots\chi_r})\nonumber\\
&=&\displaystyle\mathop{\displaystyle\mathop\sum_{\chi_1\bmod p_1}\dots
\displaystyle\mathop\sum_{\chi_k\bmod p_k}}_{\atop \chi_1\dots\chi_k\neq \chi_0}|\tau(\chi_1)|^m\dots|\tau(\chi_k)|^m\times\nonumber\\
& &\times\displaystyle\mathop{{\sum}^{\ast}}_{\chi_{k+1}\bmod p_{k+1}^{\alpha_{k+1}}}\dots\displaystyle\mathop{{\sum}^{\ast}}_{\chi_r\bmod p_r^{\alpha_r}}|\tau(\chi_{k+1})|^m\dots|\tau(\chi_r)|^m\times\nonumber\\
& &\times\quad\sum_{n=1}^{\infty}\frac{\chi_1\dots \chi_r(n)}{n(n+a)}L(1,\overline{\chi_1\dots\chi_r})\nonumber\\
&=&\displaystyle\mathop{\displaystyle\mathop\sum_{\chi_1\bmod p_1}\dots
\displaystyle\mathop\sum_{\chi_k\bmod p_k}}_{\atop \chi_1\dots\chi_k\neq \chi_0}|\tau(\chi_1)|^m\dots|\tau(\chi_k)|^m\times\nonumber\\
& &\times\displaystyle\mathop{{\sum}^{\ast}}_{\chi_{k+1}\bmod p_{k+1}^{\alpha_{k+1}}}\dots\displaystyle\mathop{{\sum}^{\ast}}_{\chi_r\bmod p_r^{\alpha_r}}\left(\sqrt{p_{k+1}^{\alpha_{k+1}}}\right)^m\dots\left(\sqrt{p_r^{\alpha_r}}\right)^m\times\nonumber\\
& &\times\quad\sum_{n=1}^{\infty}\frac{\chi_1\dots \chi_r(n)}{n(n+a)}L(1,\overline{\chi_1\dots\chi_r})\nonumber\\
&=&A\displaystyle\mathop{{\sum}^{\ast}}_{\chi_N\bmod N}N^{\frac{m}{2}}\sum_{n=1}^{\infty}\frac{\chi_1\dots \chi_r(n)}{n(n+a)}L(1,\overline{\chi_1\dots\chi_r}),
\end{eqnarray}
where
\begin{eqnarray}
A&=&\displaystyle\mathop{\displaystyle\mathop\sum_{\chi_1\bmod p_1}\dots
\displaystyle\mathop\sum_{\chi_k\bmod p_k}}_{\atop \chi_1\dots\chi_k\neq \chi_0}|\tau(\chi_1)|^m\dots|\tau(\chi_k)|^m\nonumber\\
&=&\displaystyle\mathop\sum_{\chi_1\bmod p_1}\dots
\displaystyle\mathop\sum_{\chi_k\bmod p_k}|\tau(\chi_1)|^m\dots|\tau(\chi_k)|^m-|\tau(\chi_1^0)|^m\dots|\tau(\chi_k^0)|^m\nonumber\\
&=&\left(\displaystyle\mathop{{\sum}^{\ast}}_{\chi_1\bmod p_1}|\tau(\chi_1)|^m+|\tau(\chi_1^0)|^m\right)\dots\left(\displaystyle\mathop{{\sum}^{\ast}}_{\chi_k\bmod p_k}|\tau(\chi_k)|^m+|\tau(\chi_k^0)|^m\right)\nonumber\\
& &-\left|\tau(\chi_1^0)\right|^m\dots\left|\tau(\chi_k^0)\right|^m\nonumber\\
&=&\displaystyle\mathop{{\sum}^{\ast}}_{\chi_1\bmod p_1}\displaystyle\mathop{{\sum}^{\ast}}_{\chi_2\bmod p_2}\dots\displaystyle\mathop{{\sum}^{\ast}}_{\chi_k\bmod p_k}|\tau(\chi_1)|^m\dots|\tau(\chi_k)|^m+\nonumber\\
& &+|\tau(\chi_1^0)|^m\displaystyle\mathop{{\sum}^{\ast}}_{\chi_2\bmod p_2}\dots\displaystyle\mathop{{\sum}^{\ast}}_{\chi_k\bmod p_k}|\tau(\chi_2)|^m\dots|\tau(\chi_k)|^m+\dots+\nonumber\\
& &+|\tau(\chi_k^0)|^m\displaystyle\mathop{{\sum}^{\ast}}_{\chi_1\bmod p_1}
\dots\displaystyle\mathop{{\sum}^{\ast}}_{\chi_{k-1}\bmod p_{k-1}}
|\tau(\chi_1)|^m\dots|\tau(\chi_{k-1})|^m+\nonumber\\
& &+|\tau(\chi_1^0)|^m|\tau(\chi_2^0)|^m
\displaystyle\mathop{{\sum}^{\ast}}_{\chi_3\bmod p_3}\dots\displaystyle\mathop{{\sum}^{\ast}}_{\chi_k\bmod p_k}|\tau(\chi_3)|^m\dots|\tau(\chi_k)|^m+\dots+\nonumber\\
& &+(|\tau(\chi_2^0)\dots\tau(\chi_k^0)|)^m
\displaystyle\mathop{{\sum}^{\ast}}_{\chi_1\bmod p_1}|\tau(\chi_1)|^m\nonumber\\
&=&\sum_{d|p_1\dots p_k}\displaystyle\mathop{{\sum}^{\ast}}_{\chi_M\bmod d}
|\tau(\chi_M)|^m,
\end{eqnarray}
Then we take (12) into (11) and get
\begin{eqnarray}
C_2&=&\sum_{d|M}\displaystyle\mathop{{\sum}^{\ast}}_{\chi_M\bmod d}
d^{\frac{m}{2}}N^{\frac{m}{2}}\displaystyle\mathop{{\sum}^{\ast}}_{\chi_N\bmod N}\sum_{n=1}^{\infty}\frac{\chi_M\chi_N\chi_{M}^0(n)}{n(n+a)}L(1,\overline{\chi_M\chi_N\chi_{M}^0})\nonumber\\
&=&\sum_{d|M}d^{\frac{m}{2}}N^{\frac{m}{2}}
\displaystyle\mathop{{\sum}^{\ast}}_{\chi_M\bmod d}
\displaystyle\mathop{{\sum}^{\ast}}_{\chi_N\bmod N}\sum_{n=1}^{\infty}\frac{\chi_M\chi_N\chi_{M}^0(n)}{n(n+a)}L(1,\overline{\chi_M\chi_N\chi_{M}^0})\nonumber\\
&=&\sum_{d|M}d^{\frac{m}{2}}N^{\frac{m}{2}}
\displaystyle\mathop{{\sum}^{\ast}}_{\chi \bmod Nd}
\sum_{n=1}^{\infty}\frac{\chi{\chi_{M}^0}(n)}{n(n+a)}L(1,\overline{\chi{\chi_{M}^0}}).\nonumber
\end{eqnarray}
According to Lemma 6, we easily get
\begin{eqnarray}
C_2&=&\sum_{\chi\ne\chi_0}|\tau(\chi)|^{m}\sum_{n=1}^{\infty}\frac{\chi(n)}{n(n+a)}L(1,\bar{\chi})\nonumber\\
&=&\sum_{d|M}d^{\frac{m}{2}}N^{\frac{m}{2}}\frac{\phi^{2}(N)}{N}J(d)\zeta(2)\sum_{k|q}\frac{\mu(k)}{ak^2}-\nonumber\\
&&-\sum_{d|M}d^{\frac{m}{2}}N^{\frac{m}{2}}\frac{\phi^{2}(N)}{N}J(d)\sum_{k|q}\frac{\mu(k)}{a^2k}\sum_{h=1}^{[a/k]}
\frac{1}{h}+O\left(\sum_{d|M}d^{\frac{m}{2}}N^{\frac{m}{2}}(Nd)^{\epsilon}\right)\nonumber\\
&=&N^{\frac{m}{2}-1}\phi^{2}(N)\zeta(2)\sum_{k|q}\frac{\mu(k)}{ak^2}\times\sum_{d|M}d^{\frac{m}{2}}J(d)-\nonumber\\
&&-N^{\frac{m}{2}-1}\phi^{2}(N)\sum_{k|q}\frac{\mu(k)}{a^2k}\sum_{h=1}^{[a/k]}
\frac{1}{h}\times\sum_{d|M}d^{\frac{m}{2}}J(d)+O\left(q^{\frac{m}{2}+\epsilon}\right)\nonumber\\
&=&N^{\frac{m}{2}-1}\phi^{2}(N)\zeta(2)\sum_{k|q}\frac{\mu(k)}{ak^2}
\prod_{p|M}\left(p^{\frac{m}{2}+1}-2p^{\frac{m}{2}}+1\right)-\nonumber\\
&&-N^{\frac{m}{2}-1}\phi^{2}(N)\sum_{k|q}\frac{\mu(k)}{a^2k}\sum_{h=1}^{[a/k]}
\frac{1}{h}\prod_{p|M}\left(p^{\frac{m}{2}+1}-2p^{\frac{m}{2}}+1\right)+O\left(q^{\frac{m}{2}+\epsilon}\right).\nonumber
\end{eqnarray}

(ii) The asymptotic formula of $C_3$

Following the similar method in part (i), we could get the
asymptotic formula of $C_3$, that is
\begin{eqnarray}
C_3&=&\sum_{\chi\ne\chi_0}|\tau(\chi)|^{m}
\sum_{n=1}^{\infty}\frac{\bar{\chi}(n)}{n(n+a)}L(1,\chi)\nonumber\\
&=&N^{\frac{m}{2}-1}\phi^{2}(N)\zeta(2)\sum_{k|q}\frac{\mu(k)}{ak^2}
\prod_{p|M}\left(p^{\frac{m}{2}+1}-2p^{\frac{m}{2}}+1\right)-\nonumber\\
&&-N^{\frac{m}{2}-1}\phi^{2}(N)\sum_{k|q}\frac{\mu(k)}{a^2k}\sum_{h=1}^{[a/k]}
\frac{1}{h}\prod_{p|M}\left(p^{\frac{m}{2}+1}-2p^{\frac{m}{2}}+1\right)
+O\left(q^{\frac{m}{2}+\epsilon}\right).\nonumber
\end{eqnarray}

(iii) The asymptotic formula of $C_4$

Lastly we will derive the asymptotic formula of $C_4$. Following the similar method of $C_2$ and Lemma 7, we have
\begin{eqnarray}
C_4&=&\sum_{\chi\ne\chi_0}|\tau(\chi)|^{m}
\left|\sum_{n=1}^{\infty}\frac{\chi(n)}{n(n+a)}\right|^2\nonumber\\
&=&\sum_{d|M}d^{\frac{m}{2}}N^{\frac{m}{2}}
\displaystyle\mathop{{\sum}^{\ast}}_{\chi \bmod Nd}
\left|\sum_{n=1}^{\infty}\frac{\chi{\chi_{M}^0}(n)}{n(n+a)}\right|^2\nonumber\\
&=&\sum_{d|M}d^{\frac{m}{2}}N^{\frac{m}{2}}\frac{\phi^{2}(N)}{N}J(d)\zeta(2)\sum_{k|q}\frac{\mu(k)}{a^2k^2}+\sum_{d|M}d^{\frac{m}{2}}N^{\frac{m}{2}}\frac{\phi^{2}(N)}{N}J(d)\zeta(2,\frac{a}{k})\sum_{k|q}\frac{\mu(k)}{a^2k^2}-\nonumber\\
&&-2\sum_{d|M}d^{\frac{m}{2}}N^{\frac{m}{2}}\frac{\phi^{2}(N)}{N}J(d)\sum_{k|q}\frac{\mu(k)}{a^3k}
\sum_{h=1}^{[a/k]}\frac{1}{h}+O\left(\sum_{d|M}d^{\frac{m}{2}}N^{\frac{m}{2}}(Nd)^{\epsilon}\right)\nonumber\\
&=&N^{\frac{m}{2}-1}\phi^{2}(N)\zeta(2)\sum_{k|q}\frac{\mu(k)}{a^2k^2}\times\sum_{d|M}d^{\frac{m}{2}}J(d)+\nonumber\\
&&+N^{\frac{m}{2}-1}\phi^{2}(N)\zeta(2,\frac{a}{k})\sum_{k|q}\frac{\mu(k)}{a^2k^2}
\times\sum_{d|M}d^{\frac{m}{2}}J(d)-\nonumber\\
&&-2N^{\frac{m}{2}-1}\phi^{2}(N)\sum_{k|q}\frac{\mu(k)}{a^3k}\sum_{h=1}^{[a/k]}\frac{1}{h}
\times\sum_{d|M}d^{\frac{m}{2}}J(d)+O\left(q^{\frac{m}{2}+\epsilon}\right)\nonumber\\
&=&N^{\frac{m}{2}-1}\phi^{2}(N)\zeta(2)\sum_{k|q}\frac{\mu(k)}{a^2k^2}\prod_{p|M}\left(p^{\frac{m}{2}+1}-2p^{\frac{m}{2}}+1\right)+\nonumber\\
&&+N^{\frac{m}{2}-1}\phi^{2}(N)\zeta(2,\frac{a}{k})\sum_{k|q}\frac{\mu(k)}{a^2k^2}
\prod_{p|M}\left(p^{\frac{m}{2}+1}-2p^{\frac{m}{2}}+1\right)-\nonumber\\
&&-2N^{\frac{m}{2}-1}\phi^{2}(N)\sum_{k|q}\frac{\mu(k)}{a^3k}\sum_{h=1}^{[a/k]}\frac{1}{h}
\prod_{p|M}\left(p^{\frac{m}{2}+1}-2p^{\frac{m}{2}}+1\right)+O\left(q^{\frac{m}{2}+\epsilon}\right).\nonumber
\end{eqnarray}

Combining the the estimates of (i), (ii), (iii) and Lemma 5, we get
\begin{eqnarray}
&&\sum_{\chi\ne\chi_0}|\tau(\chi)|^{m}|L(1,\chi,a)|^{2}\nonumber\\
&=&N^{\frac{m}{2}-1}\phi^{2}(N)\zeta(2,\frac{a}{k})\prod_{p||q}(1-\frac{1}{p^2})\prod_{p|M}\left(p^{\frac{m}{2}+1}
-2p^{\frac{m}{2}}+1\right)+O\left(q^{\frac{m}{2}+\epsilon}\right).\nonumber
\end{eqnarray}
This completes the proof of Theorem.\\


\begin{thebibliography}{99}
\bibitem{T.M}
M. A. Tom, Introduction to Analytic Number Theory,
Springer-Verlag, New York, 1976.
\bibitem{C.D}
C. D. Pan and C. B. Pan, Element of the Analytic Number Theory,
Science Press, Beijing, 1991 (in Chinese).
\bibitem{NA}
Ireland K and Rosen M. A Classical Introduction to Modern Number Theory. Springer-Verlag, New York, 1982, 88-91.
\bibitem{L.K}
B. C. Berndt, Generalized Dirichlet series and Hecke's functional
equation, Proc.Edinburgh Math. Soc. 15(2), 1967, 309-313.
\bibitem{S.H}
B. C. Berndt, Identities involving the coefficients of a class of
Dirichlet series. III, Trans.Amer. Math. Soc. 146, 1969, 323-342.
\bibitem{S.H}
B. C. Berndt, Identities involving the coefficients of a class of
Dirichlet series. IV, Trans. Amer. Math. Soc. 149, 1970, 179-185.
\bibitem{S.H}
R. Ma, Y. Yi, Y. L. Zhang, On the mean value of the gerneralized
Dirichlet $L$-functions, Czechoslovak Mathematical Journal 60(135)
, 2010, 597-620.
\bibitem{S.H}
R. Ma, Y. N. Niu, Y. L. Zhang, On asymptotic properties of the generalized Dirichlet $L$-functions, International Journal of Number Theory 15(6)
, 2019, 1305-1321.
\bibitem{S.H}
Y. Yi and W. P. Zhang, On the $2k$-th Power Mean of Inversion of
$L$-Functions with the Weight of Gauss Sum, Acta Mathematica Sinica, English Series 20(1), 2004, 175-180.
\bibitem{S.H}
Y. Yi and W. P. Zhang, On the $2k$-th Power Mean of Dirichlet
$L$-Functions With the Weight of Gauss Sums, Advances in
Mathematics 31(6), 2002, 517-526.
\end{thebibliography}
\end{document}